\documentclass[%
 reprint,
nofootinbib,
 amsmath,amssymb,
 aps, prl
]{revtex4-1}
\usepackage{graphicx}
\usepackage{dcolumn}
\usepackage{bm}
\usepackage{diagbox}
\usepackage{subcaption}
\usepackage{pifont}

\begin{document}

\preprint{APS/123-QED}

\title{Causal Stability and Synchronization}

\author{Aditi Kathpalia}
\email{kathpaliaadit@nias.res.in}
\author{Nithin Nagaraj}%
 \email{nithin@nias.res.in}
\affiliation{%
Consciousness Studies Programme \\
National Institute of Advanced Studies\\
Indian Institute of Science Campus, Bengaluru, India
}%




\date{\today}

\begin{abstract}
Synchronization of chaos arises between coupled dynamical systems and is very well understood as a temporal phenomena which leads the coupled systems to converge or develop a dependence with time. In this work, we provide a complementary spatial perspective to this phenomenon by introducing the novel idea of causal stability. We then propose and prove a causal stability synchronization theorem as a necessary and sufficient condition for synchronization. We also provide an empirical criteria to identify synchronizing variables in coupled identical chaotic dynamical systems based on causality analysis on time series data of the driving system alone. \\
%
%

\end{abstract}                            
\maketitle


\section{\label{sec:level1}Introduction}

Synchronization of coupled chaotic systems is a well-known and well-studied phenomenon~\cite{pikovsky2003synchronization}. Pecora and Carroll's 1990 paper~\cite{pecora1990synchronization} followed by He and Vaidya's work in 1992~\cite{he1992analysis} were a revelation into chaotic synchronization and opened up an entire field of intense research. Chaotic synchronization has found applications in living systems~\cite{glass1988clocks, mosekilde2002chaotic}, human cognition and neuroscience~\cite{buzsaki2006rhythms, rodriguez1999perception, singer2011consciousness} as well as physics, chemistry and engineering~\cite{strogatz1994nonlinear}. Eventually a lot of work followed up on tracking the transition from incoherence to synchrony~\cite{rosenblum1996phase, boccaletti2002nonlinear, pikovsky1997phase, paluvs2001synchronization, lahav2018synchronization}.

Though Pecora and Carroll, He and Vaidya proved temporal conditions for synchronization (negative conditional Lyapunov exponents and asymptotic stability of the non-driven subsystem), there exists no spatial perspective of the same. For ergodic dynamical processes, it is only natural to expect a spatial counterpart for temporal conditions of chaotic synchronization. In this work, we explore causal interactions within the master and slave systems to identify spatial influences that drive a slave system to synchronize with its master.

Sensitive dependence on initial conditions, the hallmark of chaos, implies that trajectories from nearby initial conditions diverge.  However, when these chaotic systems are coupled, the \emph{causal influence} of the master on the slave via information transmitted by the coupled variable may lead to synchronized behaviour. The causal influence of the coupled variable on the non-driven subsystem is of importance in this regard. If this influence is invariant to perturbation of the subsystems' initial conditions then we can expect chaotic synchronization.

Asymptotic stability~\cite{he1992analysis} - the condition that the non-driven subsystem reaches the same eventual state at a fixed time no matter what the initial conditions were, is both a necessary and sufficient for chaotic synchronization. We are interested in asking \emph{what is the cause of asymptotic stability of the subsystem?} Specifically, we explore what kind of causal influence does the forced variable have on this subsystem that leads it to be driven to the same eventual state each time. We formally derive and prove a necessary and sufficient condition for chaotic synchronization based on causal influence to the non-driven subsystem. 

Furthermore, for identical master and slave systems, given the dynamics of the master alone,  we provide an empirical criteria for determining which subsystem will be driven to synchronization, or in other words, which variables when coupled will result in synchronization. This is done based on intra-system causal influences estimated entirely from time series data of the master system. This is an important novel contribution with real world applications in the control of chaos, especially when we need to determine which nodes should be coupled for synchronization when the underlying equations are completely unknown.

An important tool for the above discussed causal analysis is a recently proposed measure `Compression-Complexity Causality (CCC)'~\cite{kathpalia2019data} which efficiently captures causal relationships between time series of coupled processes based on dynamical complexity. The measure provides not only the \emph{quantity} (strength) of causality but also its \emph{quality} which is reflected in the sign of the estimated $CCC$ value. This is an important property that captures the \emph{kind of dynamical influence} from one variable to another. The sign of $CCC$ determines whether the influence to the driven variable from its own past is different from or similar to that from the past of the driving variable. The similar or dissimilar dynamics has the potential to determine whether the slave will be driven to synchronize with the master or will remain `sensitive to its own initial conditions', which are different from that of the master.


\section{Synchronization via Causal Stability}

Let a master system be governed by the following set of differential equations:
\renewcommand{\vec}[1]{\mathbf{#1}}
\begin{equation}
\label{eq_master_tog}
    \dot{\vec{y}}=\vec{f}(t,\vec{y}),
\end{equation}
where $\dot{\vec{y}}$ and $\vec{f}$ are vectors. As described in~\cite{he1992analysis}, the master system (Eq.~\ref{eq_master_tog}) can be divided into two interdependent subsystems $\vec{p}$ and $\vec{q}$. The slave system also consists of two subsystems $\vec{p'}$ and $\vec{q'}$, whose functional form is identical to the corresponding master system. $\vec{p'}$ component of slave is completely overridden by the $\vec{p}$ component of master and will be referred to as the forced/driven variable. $\vec{q'}$ component is allowed to have initially different conditions and will be referred to as the non-driven subsystem. The equations for the master are as below:
\begin{equation}
    \dot{\vec{p}}=\vec{h}(t,\vec{p},\vec{q}), ~~~~ \dot{\vec{q}}=\vec{g}(t,\vec{p},\vec{q}), 
    \label{eq_master}
\end{equation}
and the slave system is given by:
\begin{equation}
    \vec{p'}=\vec{p}, ~~~~ \dot{\vec{q'}}=\vec{g}(t,\vec{p},  \vec{q'}).
    \label{eq_slave}
\end{equation}
As pointed before, the influence of the forced variable ($\vec{p'}$) on the non-driven subsystem ($\vec{q'}$) is important for our analysis of the behavior of slave dynamics. `A \textit{cause} or \textit{causal influence} is defined as something that makes a difference and the difference it makes must be a difference from what would have happened without it'. This definition was given by philosopher David Lewis inspired by David Hume's notion of causality~\cite{pearl2018book}. For evaluation of cause or causal influence based on time series analysis, the above definition of causality was formulated in a mathematical way by Wiener~\cite{wiener}. According to Wiener, a time series $Y$ is a cause for time series $X$ if the past of $Y$ contains information that helps predict $X$ above and beyond the information contained in past values of $X$ alone. This principle led to Wiener-Granger causality measure~\cite{granger} for coupled autoregressive processes. Subsequently, concepts related to information flow, such as Transfer Entropy~\cite{schreiber} or Conditional Mutual Information~\cite{paluvs2001synchronization} were developed for causality testing. Recently, we have proposed a measure called Compression-Complexity Causality (CCC)~\cite{kathpalia2019data} that captures causality based on dynamical evolution of processes.

%

Let us suppose the subsystem $\vec{p}$ comprises of the variable $x$ and $\vec{p'}$ of variable $x'$. Also, the subsystem $\vec{q}$ comprises of the variables $y$ and $z$ and $\vec{q'}$ of variables $y'$ and $z'$. We estimate the net causal influence input to the subsystem $\vec{q'}$ as below. The causal influence can be determined using any of the methods discussed above and is denoted by $CI$.
%
%
\begin{equation}
    CI_{net(\vec{p'} \rightarrow \vec{q'})}=CI_{\vec{p'} \rightarrow \vec{q'}} - CI_{\vec{q'} \rightarrow \vec{p'}},
\end{equation}
which in this case reduces to:
\begin{equation}
\begin{split}
    CI_{net(x' \rightarrow y',z')}=CI_{x' \rightarrow y'|z'} + CI_{x' \rightarrow z'|y'} \\- CI_{y' \rightarrow x'|z'} - CI_{z' \rightarrow x'|y'}.
    \end{split}
\end{equation}
where $CI_{a \rightarrow b|c}$ is the conditional causal influence of $a$ on $b$ given time series $c$. In case of more than three variables, the conditioning is performed on all the remaining variables (other than $a$ and $b$).

We define \textit{causal stability} for the subsystem $\vec{q'}$ as follows: 

\noindent \textbf{Definition 1} {\it Causal Stability:} The subsystem $\vec{q'}$ is said to be causally stable if:
\begin{equation}
|CI_{net(\vec{p'}(\vec{p'_0}) \rightarrow \vec{q'}(\vec{q'_0}))}-CI_{net(\vec{p'}(\vec{p'_0}) \rightarrow \vec{q'}(\vec{q'_1}))}| < \epsilon,
\label{eq_cs}
\end{equation}
for arbitrarily small $\epsilon>0$.  $\vec{p'}(\vec{p'_0})$ represents the solution trajectory of system $\vec{p'}$ started from initial vector $\vec{p'_0}$, that is $\vec{p'(t_1:t_2;t_0,p'_0)}$, for all $t_1$ such that $t_{tr}<t_1<+\infty$, $t_2=t_1+T$, $T>0$ and transients for the trajectory last until $t_{tr}$. $T$ is the length of the trajectory over which $CI_{net}$ is estimated. $\vec{q'}(\vec{q'_0})$ and $\vec{q'}(\vec{q'_1})$ represent the solution trajectory of system $\vec{q'}$ started from initial vectors $\vec{q'_0}$ and $\vec{q'_1}$ respectively, that is $\vec{q'(t_1:t_2;t_0,q'_0)}$ and $\vec{q'(t_1:t_2;t_0,q'_1)}$, for the same $t_1, t_2$ defined above.

If the solution trajectory for $\vec{q'}$ in Eq.~\ref{eq_slave} is causally stable  $\forall~\vec{q'_0} \in D(t_0)$, where $D(t_0) \subseteq R^n$, then $D(t_0)$ is called the region of causal stability. If $D(t_0)=R^n$, then the solutions are said to be \emph{globally causally stable}. The requirement on $D(t_0)$ is that any two initial conditions for the slave subsystem $\vec{q'}$ taken from this region yield solution trajectories which at some finite point in time ($\geq t_{tr}$) come close to each other for some chunk of time. Intuitively, causal stability implies that the net causal influence input to the non-driven subsystem $\vec{q'}$ from the driven subsystem $\vec{p'}$ is \emph{invariant} to the change in the initial conditions of $\vec{q'}$.

We define synchronization as:\\
\noindent \textbf{Definition 2} \textit{Synchronization}: Let us take two systems $\dot{\vec{y}}=\vec{f}(t,\vec{y})$ and $\dot{\vec{y'}}=\vec{f'}(t,\vec{y'})$, where $\vec{y}, \vec{y'} \in R^n$. Let their solutions be given by $\vec{y}(t;t_0,\vec{y_0})$ and $\vec{y'}(t;t_0,\vec{y'_0})$, respectively. $\vec{f}(t,\vec{y})$ synchronizes with $\vec{f'}(t,\vec{y'})$ if there exists $D(t_0) \subseteq R^n$, such that $\vec{y_0}, \vec{y'_0} \in D(t_0)$ implies
\begin{equation}
    \|\vec{y}(t;t_0,\vec{y_0})-\vec{y'}(t;t_0,\vec{y'_0})\| \rightarrow 0 \mathtt{~as~} t \rightarrow \infty.
\end{equation}
We state and prove the following theorem:
\noindent \par \textbf{Causal Stability Synchronization Theorem:} The slave system $(\vec{p'}, \vec{q'})$ synchronizes with the master system $(\vec{p}, \vec{q})$ \emph{iff} there exists $D(t_0) \subseteq R^n$ such that when the initial conditions of the non-driven part of the slave system $\dot{\vec{q'}}=\vec{g}(t,\vec{p},\vec{q'})$ fall in $D(t_0)$, the solution trajectory of $\vec{q'}$ is causally stable. 

\noindent \textit{Proof}:

\textit{Sufficient Condition}: If there exists a $D(t_0) \subseteq R^n$ such that when the initial conditions of the non-driven part of the slave system $\dot{\vec{q'}}=\vec{g}(t,\vec{p},\vec{q'})$ fall in $D(t_0)$, the solution trajectory of $\vec{q'}$ is causally stable, then, by the definition of causal stability,
\begin{equation}
|CI_{net(\vec{p'}(\vec{p'_0}) \rightarrow \vec{q'}(\vec{q'_0}))}-CI_{net(\vec{p'}(\vec{p'_0}) \rightarrow \vec{q'}(\vec{q'_1}))}| < \epsilon,
\label{eq_cs2}
\end{equation}
where $\vec{q'_0}$ and $\vec{q'_1}$ are two arbitrary initial conditions in $D(t_0)$. 
Let us suppose the net difference in causal influence that the presence of a variable $\vec{p'}$ makes on the future of $\vec{q'}$ above and beyond the past trajectory of $\vec{q'}$ is computed using a measure  $C$~\footnote{$C$ could be an infotheoretic quantity such as conditional entropy or could be a complexity measure.
}. Thus the Eq.~\ref{eq_cs2} can be elaborated as below:
\begin{equation}
\label{eq_cs_diff}
{\small
\begin{split}
    C(\Delta \vec{q'}(\vec{q'_0})| \vec{q'}(\vec{q'_0})_{past}, \vec{p'}(\vec{p'_0})_{past}) - C(\Delta \vec{q'}(\vec{q'_0})| \vec{q'}(\vec{q'_0})_{past}) \\
    \approx C(\Delta \vec{q'}(\vec{q'_1})| \vec{q'}(\vec{q'_1})_{past}, \vec{p'}(\vec{p'_0})_{past}) - C(\Delta \vec{q'}(\vec{q'_1})| \vec{q'}(\vec{q'_1})_{past}),
    \end{split}
    }
\end{equation}

where $\Delta \vec{q'}(\vec{q'_0})$ and $\Delta \vec{q'}(\vec{q'_1})$ are the current window of time series data from the system $\vec{q'}$ started from initial conditions $\vec{q'_0}$ and $\vec{q'_1}$ respectively. $\vec{q'}(\vec{q'_0})_{past}$ and  $\vec{q'}(\vec{q'_1})_{past}$ represent the immediate past values of the window from the system $\vec{q'}$ started from initial conditions $\vec{q'_0}$ and $\vec{q'_1}$ respectively. $\vec{p'}(\vec{p'_0})_{past}$ represent synchronous past values from the system $\vec{p'}$ started from initial condition $\vec{p'_0}$.

If the system $\vec{q'}$ is started from two nearby initial conditions, then their initial solution trajectories will be similar as the influence from the fixed master $\vec{p}(\vec{p_0})$ ($= \vec{p'}(\vec{p'_0})$) time series has not yet come into play. Even if the initial conditions are far away, then the condition on $D(t_0)$ would require the trajectories to come close to each other for a chunk of time starting at some finite time $t_s$. Then, for this chunk of time, the trajectories are close to each other and hence:
\begin{equation}
    C(\vec{q'}(\vec{q'_0})_{past})_{t_s}=C(\vec{q'}(\vec{q'_1})_{past})_{t_s}.
    \label{eq_cs_cond1}
\end{equation}
Also,
\begin{equation}
    C(\Delta \vec{q'}(\vec{q'_0})|\vec{q'}(\vec{q'_0})_{past})_{t_s}=C(\Delta \vec{q'}(\vec{q'_1})|\vec{q'}(\vec{q'_1})_{past})_{t_s}.
    \label{eq_cs_cond2}
\end{equation}
The reason for the above two equalities for measure $C$ is that it is computed using coarse grained time series (symbolic sequences) or k-nearest neighbor estimation techniques, rendering the measure to become equal for two close by trajectories. Eqs.~\ref{eq_cs_cond1} and \ref{eq_cs_cond2} hold at time $t_s$ and Eq.~\ref{eq_cs_diff} (and~\ref{eq_cs2}) hold for all $t \geq t_s$. So, Eq.~\ref{eq_cs_diff} is true at $t_1=t_s+ \Delta t$ where $\Delta t$ is an arbitrary increment in time. For this to be true, conditions~\ref{eq_cs_cond1} and \ref{eq_cs_cond2} need to also hold at $t_1$, because in order to maintain~\ref{eq_cs2} for every small time step increment $\Delta t$ starting at $t_s$, the measure $C$ for the solution trajectories started at two different initial conditions cannot drastically change, nor can the influence in each case from its own past to future. Also, the only drivers of system $\vec{q'}$, or specifically the contributors to the evolution of $\vec{\Delta q'}$ are the past of $\vec{q'}$ itself and the forced variable $\vec{p'}$ (which is fixed across different initial conditions of $\vec{q'}$). Since the measure $C$ for the solution trajectories $\vec{q'}(\vec{q'_0})$ and $\vec{q'}(\vec{q'_1})$ is equal as well as the influence they bring to their future, these solution trajectories have been forced not to diverge but to converge as they were similar at $t_s$ and need to maintain these conditions at $t_1$. Thus,
%
%
\begin{equation}
    \|\vec{q'}(t,\vec{p};t_0,\vec{p_0},\vec{q'_0})-\vec{q'}(t,\vec{p};t_0,\vec{p_0},\vec{q'_1})\| \rightarrow 0 \mathtt{~as~} t\rightarrow \infty.
\end{equation}
%





\textit{Necessary Condition}: If there exists a $D(t_0) \subseteq R^n$ such that when the initial conditions of the non-driven part of the slave system $\dot{\vec{q'}}=\vec{g}(t,\vec{p},\vec{q'})$ fall in $D(t_0)$, the master and slave systems are synchronized. Suppose the master starts from initial condition $\vec{q_m}$ and $\vec{q'_0}$ and $\vec{q'_1}$ are two arbitrary initial conditions taken for the slave in the set $D(t_0)$. Then, by the definition of synchronization stated above,
\begin{equation}
    \|\vec{q'}(t,\vec{p};t_0,\vec{p_0},\vec{q'_0})-\vec{q}(t,\vec{p};t_0,\vec{p_0},\vec{q_m})\| \rightarrow 0, t\rightarrow \infty,
    \label{eq_synch_proof1}
\end{equation}
\begin{equation}
    \|\vec{q'}(t,\vec{p};t_0,\vec{p_0},\vec{q'_1})-\vec{q}(t,\vec{p};t_0,\vec{p_0},\vec{q_m})\| \rightarrow 0, t\rightarrow \infty.
    \label{eq_synch_proof2}
\end{equation}

Using Eqs.~\ref{eq_synch_proof1} and~\ref{eq_synch_proof2},
\begin{equation}
    \|\vec{q'}(t,\vec{p};t_0,\vec{p_0},\vec{q'_0})-\vec{q'}(t,\vec{p};t_0,\vec{p_0},\vec{q'_1})\| \rightarrow 0, t\rightarrow \infty.
\end{equation}
Because of the uniqueness of solutions,
\begin{equation}
    \vec{q'}(t,\vec{p};t_0,\vec{p_0},\vec{q'_0})=\vec{q'}(t,\vec{p};t_0,\vec{p_0},\vec{q'_1}).
\end{equation}
%
So, $CI_{net}$ estimated from a fixed time series (for solution trajectory starting at any time $t$) to the two identical time series of system $\vec{q'}$ will be equal. Hence,
\begin{equation}
|CI_{net(\vec{p'}(\vec{p'_0}) \rightarrow \vec{q'}(\vec{q'_0}))}-CI_{net(\vec{p'}(\vec{p'_0}) \rightarrow \vec{q'}(\vec{q'_1}))}| = 0,
\label{eq_cs3}
\end{equation}
which implies causal stability of $\vec{q'}$.$\hfill \square$
\subsection{Empirical Study}
The difference of $CI_{net}$ to the non-driven subsystem of slave for two different initial conditions approaches zero when the system is led to synchronization and not otherwise, is demonstrated for the Lorenz system in Fig.~\ref{fig_lorenz_res}, whose simulation and $CI_{net}$ estimation details are given below. $CI_{net}$ is estimated using the CCC measure. CCC happened to be the most appropriate choice as it is a model-free measure of causality applicable for non-linear time series and has been tested on several real-world like scenario simulations as well as real world datasets~\cite{kathpalia2019data,agarwal2019distinguishing}. Moreover, for our purpose we need a measure, using which, causal influence can be estimated over short lengths of windows taken from time series data. Info-theoretic measures based on probability density estimates do not perform very well in this regard. Furthermore, CCC is an interventional data-based causality measure based on dynamical evolution of processes and is not merely based on associational relations. In this regard, it is more faithful to the underlying mechanism generating the data. Lorenz system was simulated using Euler's method as per the following equations:  
\begin{equation}
\label{eq_lorenz_sys}
\begin{split}
{\frac {\mathrm {d} x}{\mathrm {d} t}}&=\sigma (y-x),\\{\frac {\mathrm {d} y}{\mathrm {d} t}}&=x(\rho -z)-y,\\{\frac {\mathrm {d} z}{\mathrm {d} t}}&=xy-\beta z,
\end{split}
\end{equation}
where, $\sigma=10$, $\rho=60$ and $\beta=8/3$.

\begin{figure*}[t]
\centering
\includegraphics[width=\textwidth]{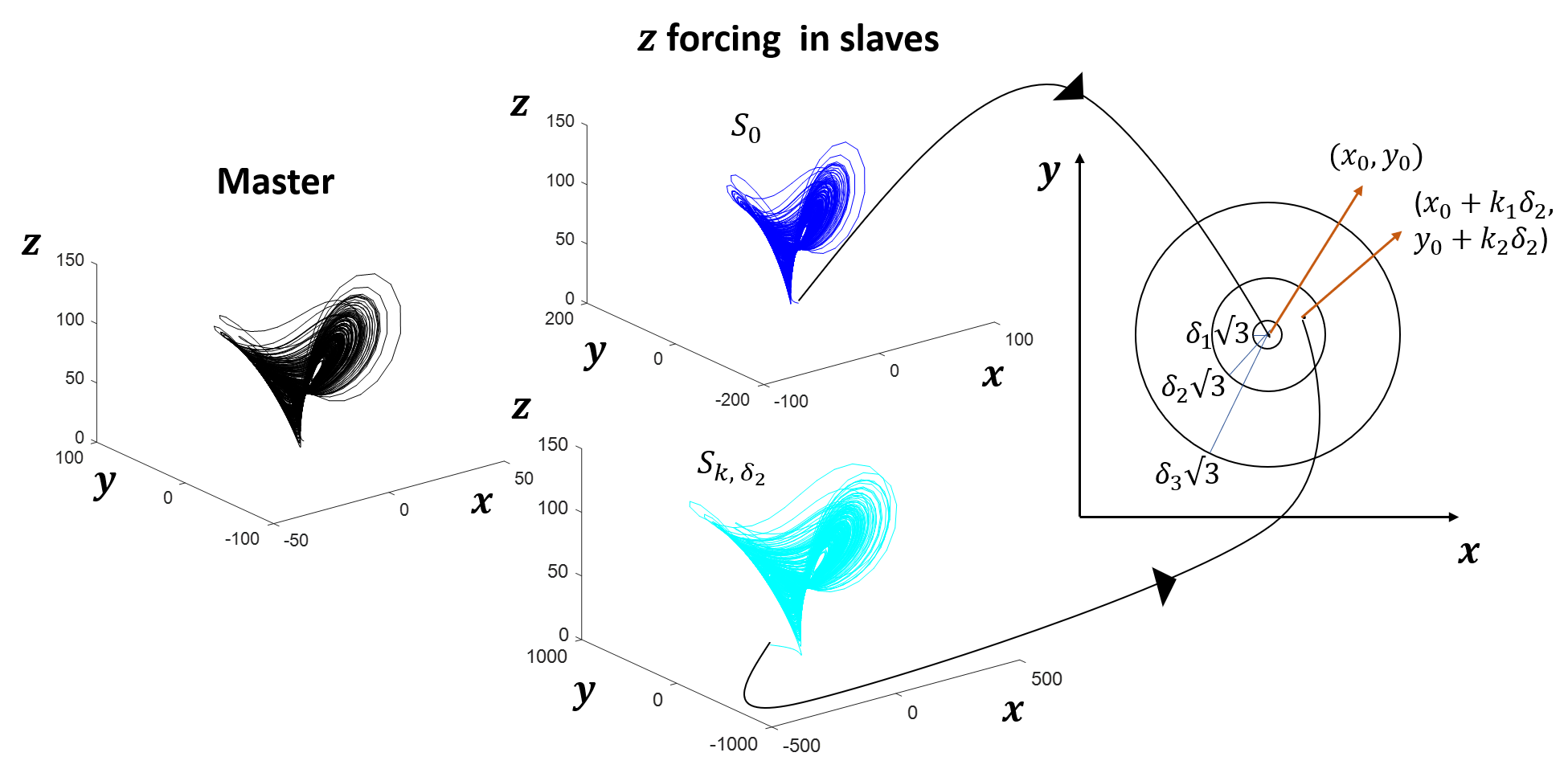}
\caption{(color online). Master Lorenz attractor (leftmost, black color). Slave Lorenz $S_0$ with forced $z$ variable and fixed initial conditions $(x_0, y_0)$ for $x$ and $y$ (middle, top, blue color). Slave Lorenz $S_{k,\delta_2}$ with forced $z$ variable and initial conditions in a $\delta_2\sqrt3$ radius disc around $(x_0, y_0)$ (middle, bottom, cyan color). Initial conditions for $x$, $y$ variables of both the slaves depicted in rightmost figure.}
\label{fig_lorenz}
\end{figure*}

Master Lorenz ($x,y,z$) was simulated starting from initial conditions $(x_m,y_m,z_m)=(3,4,6)$. For the slaves' dynamics ($x',y',z'$), one of the variables was forced to be the same as the master -- either $x=x'$ or $y=y'$ or $z=z'$. We fix one of the slaves $S_0$ to start with initial conditions: $(x'_0,y'_0,z'_0)=(7,1,6)$. The second slave $S_{k,\delta}$ was started within a sphere of radius $\delta\sqrt{3}$ from the first slave, its initial conditions are given by $(x'_{k,\delta},y'_{k,\delta},z'_{k,\delta})=(x_0+k_1\delta,y_0+k_2\delta,z_0+k_3\delta)$ where $k_1, k_2, k_3$ are independently chosen uniformly at random from the set $(-1, 1)$. As an example, for $z$-forcing in Lorenz system, the master and two slave attractors as well as their initial condition are depicted in Fig.~\ref{fig_lorenz}. 10,000 time points were simulated for both the master and the slave after removal of 2000 samples (transients).  For three different settings of $\delta$ ($\delta_1=1, \delta_2=10, \delta_3=100$), we simulate several \emph{secondary} slaves $S_{k,\delta}$ ($k=1,2, \ldots, 100$). For increasing $k$, we estimate and plot the mean of the absolute differences in the $CCC_{net}$ values between $S_0$ and secondary slaves in Fig.~\ref{fig_lorenz_res} with appropriately chosen parameters for $CCC$ estimation\footnote{The parameters used in the computation of $CCC$ in case of Lorenz are $L=150$, $w=15$, $\delta=80$, $B=8$, in case of {R{\"o}ssler} are $L=300$, $w=15$, $\delta=200$, $B=8$, in case of 5D system are $L=450$, $w=80$, $\delta=300$, $B=4$ and in case of Chen and  H{\'e}non are $L=100$, $w=15$, $\delta=80$, $B=8$. These parameters were selected on the basis of parameter selection criteria and rationale given in the supplementary material of~\cite{kathpalia2019data}.}. The mean is given by the following expression:
\begin{equation}
\label{eq_mean_CCC_net}
M_{k,\delta}=\frac{1}{k} \sum_{i=1}^{k} |CCC_{net}(S_0) - CCC_{net}(S_{i,\delta})|,
\end{equation}
where 
\begin{equation}
    \begin{split}
CCC_{net}(S_0) &= CCC_{net({z'}({z_m}) \rightarrow ({x'(x'_0),y'(y'_0)})}, \\
CCC_{net}(S_{i,\delta}) &= CCC_{net({z'}({z_m}) \rightarrow ({x'(x'_{k,\delta}),y'(y'_{k,\delta})})}, 
\end{split}
\end{equation}
when $z$ is forced ($z=z'$). The above analysis is done independently for $x$ and $y$ forcing as well.

We see that since $x$ and $y$ forcing lead to synchronization of the slaves with the master, the difference in the $CCC_{net}$ values turns out to be zero for any slave chosen from either of the discs. On the other hand, since $z$ forcing does not lead to synchronization, mean $CCC_{net}$ difference is non-zero and increases with increase in $\delta$. From Fig.~\ref{fig_lorenz_res}, we infer that $z$ forcing results in a \emph{causally unstable} non-driven slave subsystem whereas $x$ (or $y$) forcing leads to a \emph{causally stable} non-driven slave subsystem. This successfully validates \emph{causal stability synchronization theorem} for the Lorenz system.

\begin{figure*}[t]
\centering
\includegraphics[width=\textwidth]{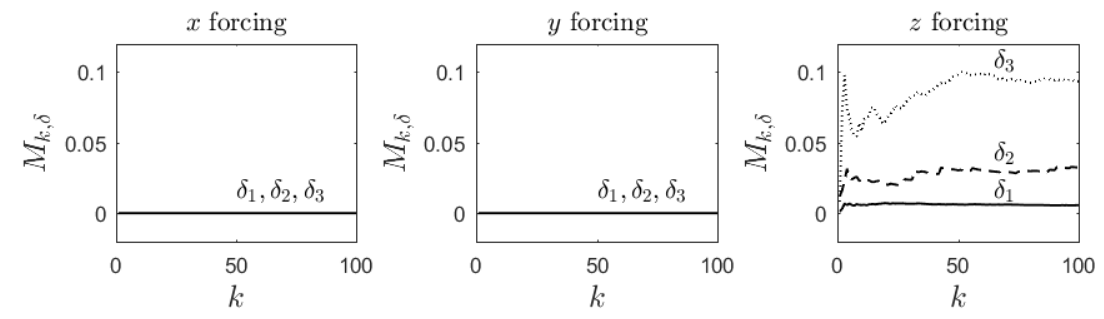}
\caption{Mean of absolute $CCC_{net}$ differences, $M_{k, \delta}$, between a fixed slave Lorenz $S_0$ and secondary slave Lorenz $S_{k,\delta}$ for ${k=1,2,\ldots,100}$ within a $\delta\sqrt3$ radius of $S_0$;  $\delta_1=1$ (solid line), $\delta_2=10$ (dashed line) and $\delta_3=100$ (dotted line). Causal Stability satisfied for $x$ and $y$ forcing but not for $z$ forcing.}
\label{fig_lorenz_res}
\end{figure*}

While Eq.~\ref{eq_cs2} is proven to be the necessary and sufficient condition for complete synchronization, we conjecture that a more relaxed causal stability condition shall be true for generalized synchronization~\cite{rulkov1995generalized, kocarev1996generalized}:
\begin{equation}
|CI_{net(\vec{p'}(\vec{p'_0}) \rightarrow \vec{q'}(\vec{q'_0}))}-CI_{net(\vec{p'}(\vec{p'_0}) \rightarrow \vec{q'}(\vec{q'_1}))}| \leq f(\overrightarrow{\delta}),
\label{eq_cs_generalized}
\end{equation}
where $\overrightarrow{\delta}$ is the magnitude of difference vector between initial conditions $\vec{q'_0}$ and $\vec{q'_1}$. The exact form of $f$ needs to be determined which is outside the scope of this paper.

\section{Synchronizing Variables}

As per the \emph{causal stability synchronization theorem}, if the net input causal influence to the non-driven subsystem is invariant with change in initial conditions, those initial conditions are led to synchronize with the master. The net input causal influence to the non-driven subsystem happens to be the net output causal influence from the driven subsystem to the non-driven part. It can be said that, this particular influence is responsible for synchronization. For a slave which is synchronized with the master, the net causal influence from the driven to non-driven subsystem in the master and slave remains the same. It should be possible then to decide on the basis of the \emph{nature of intra-system causal influences} within the master, coupling of which variables may lead to synchronization\footnote{We term these as \emph{synchronizing variables}.}. Given a network of several coupled dynamical systems, based on the time series of the driver system alone, the \textit{causality} perspective that we have proposed in this paper, can resolve an important question -- \emph{what specific properties do synchronizing variables exhibit?} This kind of analysis can be very useful for networks where we wish to control chaos by adjusting magnitude and direction of coupling between systems as well as the selection of coupling variables.

To address this important question, we use the \emph{sign of CCC} measure, since it gives information on the `kind of dynamical causal influence' which the \emph{cause} variable has on the \emph{effect} variable. If the kind of dynamical influence from a variable $x$ to another variable $y$ is different from the past of $y$ to itself, then $CCC_{x \rightarrow y} < 0$ . On the other hand, if the kind of dynamical influence from a variable $x$ to another variable $y$ is similar to that from the past of $y$ to itself, then $CCC_{x \rightarrow y} > 0$. This is true also for conditional causality estimation when more than two variables are there in the system. For further details on negative CCC please refer to~\cite{kathpalia2019data} and its supplementary material. 

While there can be various mechanisms for coupling dynamical systems to study chaotic synchronization, we consider the simplest case where forcing a particular variable in the slave system to become identical to that of the master system, may result in complete synchronization for all the variables. The slave and master are taken to be identical systems apart from their initial conditions. 

Lorenz system was simulated as per Eq.~\ref{eq_lorenz_sys} using Euler's method, where, $\sigma=10$, $\rho=60$ and $\beta=8/3$, which is known to exhibit chaos. In this case, it is well-known that forcing either $x$ or $y$ leads to complete synchronization whereas $z$ forcing does not~\cite{he1992analysis}. Table~\ref{table_5_systems}(a) shows conditional intra-system causality values between each pair of variables as well as $CCC_{net}$ values from each variable to the corresponding subsystem. We see that the net causal influence, $CCC_{net(x\rightarrow y,z)}<0$ as well as $CCC_{net(y\rightarrow x,z)}<0$ while $CCC_{net(z\rightarrow x,y)}>0$. We have an intuitive understanding for why this happens. The synchronizing variables ($x$ and $y$) influence their corresponding non-driven subsystems with a negative $CCC$ bringing a dynamical influence on the subsystem which is different from its own past. This kind of an influence constrains the subsystem of the slave to not follow its own past dynamics and is driven by the forced variable towards complete synchronization. $z$ brings an influence on the ($x$, $y$) subsystem which is commensurate with the subsystem's own past. Thus, influence from $z$ is unable to constrain $x$ and $y$ adequately to override their own past dynamics and thus is unable to result in synchronization.

This behavior was studied in a number of other continuous-time dynamical systems such as -- R{\"o}ssler, Chen and a 5D system in the chaotic regime.

\textit{R{\"o}ssler}:
\begin{equation}
\begin{split}
{\frac{dx}{dt}}=-y-z\\{\frac{dy}{dt}}=x+ay
\\{\frac{dz}{dt}}=b+z(x-c),
\end{split}
\end{equation}
where $a=0.2$, $b=0.2$ and $c=9$.

\textit{Chen}:
\begin{equation}
\begin{split}
{\frac{dx}{dt}}=a(y-x)
\\{\frac{dy}{dt}}=(c-a)x-xz+cy
\\{\frac{dz}{dt}}=xy-bz,
\end{split}
\end{equation}

where $a=35$, $b=3$, $c=28$.

where $a=1.4$ and $b= 0.3$.

\begin{equation}
\label{eq_5D_lorenz_sys}
\begin{split}
{\frac {\mathrm {d} x}{\mathrm {d} t}}&=\sigma (y-x)+w,\\{\frac {\mathrm {d} y}{\mathrm {d} t}}&=x(\rho -z)-y,\\{\frac {\mathrm {d} z}{\mathrm {d} t}}&=xy-\beta z, \\{\frac {\mathrm {d} q}{\mathrm {d} t}}&=-q^3+w, \\{\frac {\mathrm {d} w}{\mathrm {d} t}}&=-x-q-8w.
\end{split}
\end{equation}
where, $\sigma=10$, $\rho=60$ and $\beta=8/3$.

Though the theorem is proved for continuous time systems, the sign of $CCC$ values was analyzed to identify synchronizing variables even for discrete time systems. It has been shown earlier that when two one-dimensional Tent map systems are coupled, we get negative $CCC$ from the independent map to the dependent map, showing that the kind of causal influence is different from the past of the independent map~\cite{kathpalia2019data}. In order to identify synchronizing variables for discrete time systems using intra-system $CCC$ values, we simulated the well-known 2D H{\'e}non map in the chaotic regime,

\textit{H{\'e}non}:
\begin{equation}
\begin{split}
x_{n+1}=1-x_{n}^{2}+y_{n}\\y_{n+1}=bx_{n},
\end{split}
\end{equation}
where $n$ stands for discrete time. 
\begin{table}[h!]
\centering
\renewcommand{\arraystretch}{1.3}
\newcolumntype{C}[1]{>{\centering\arraybackslash}m{#1}}
\caption{Conditional $CCC$ values estimated between variables of (a) Lorenz, (b) R{\"o}ssler, (c) Chen and (d) 5D continuous-time system, (e) H{\'e}non and corresponding $CCC_{net}$ values from each variable to its subsystem.}
\label{table_5_systems}
\subcaptionbox{Lorenz}{
\begin{tabular}
{|c|C{1.5cm}|C{1.5cm}|C{1.5cm}|}\hline
\diagbox{To}{From} & {\bf $x$} & {\bf $y$} & {\bf $z$} \\ \hline
{\bf $x$} & 0 & -0.0270 & 0.0390  \\ \hline
{\bf $y$} & -0.0250 & 0 & 0.0330 \\ \hline
{\bf $z$} & 0.0251 & -0.0040 & 0 \\ \hline
{\bf $CCC_{net}$} & {\bf -0.0119} & {\bf -0.0390} & {\bf 0.0509}
\\\hline
\end{tabular}
}
\quad
\subcaptionbox{R{\"o}ssler}{
\begin{tabular}{|c|C{1.5cm}|C{1.5cm}|C{1.5cm}|}\hline
\diagbox{To}{From} & {\bf $x$} & {\bf $y$} & {\bf $z$} \\ \hline
{\bf $x$} & 0 & 0.0412 & 0.0730  \\ \hline
{\bf $y$} & 0.0345 & 0 & 0.0724 \\ \hline
{\bf $z$} & 0.0337 & 0.0288 & 0
\\\hline
{\bf $CCC_{net}$} & {\bf -0.0460} & {\bf -0.0369} & {\bf 0.0829}
\\\hline
\end{tabular}
}
\subcaptionbox{Chen}{
\begin{tabular}{|c|C{1.5cm}|C{1.5cm}|C{1.5cm}|}\hline
\diagbox{To}{From} & {\bf $x$} & {\bf $y$} & {\bf $z$} \\ \hline
{\bf $x$} & 0 & -0.0388 & 0.0769  \\ \hline
{\bf $y$} & -0.0271 & 0 & 0.0814 \\ \hline
{\bf $z$} & 0.0299 & 0.0311 & 0 \\ \hline
{\bf $CCC_{net}$} & {\bf -0.0353} & {\bf -0.0620} & {\bf 0.0973}
\\\hline
\end{tabular}
}
\subcaptionbox{5D continuous-time system}{
\begin{tabular}{|c|C{1.18cm}|C{1.18cm}|C{1.18cm}|C{1.18cm}|C{1.18cm}|}\hline
\diagbox{To}{From} & {\bf $x$} & {\bf $y$} & {\bf $z$} & {\bf $q$} & {\bf $w$} \\ \hline
{\bf $x$} & 0 & -0.0275 & 0.0163 & 0.0326 & 0.0287\\ \hline
{\bf $y$} & -0.0160 & 0 & 0.0170 & 0.0356 & 0.0318\\ \hline
{\bf $z$} & 0.0090 & -0.0021 & 0 & 0.0406 & 0.0413 \\ \hline
{\bf $q$} & 0.0020 & 0.0009 & 0.0296 & 0 & -0.0001 \\ \hline
{\bf $w$} & -0.0018 & -0.0069 & 0.0291 & 0.0098 & 0 \\ \hline
{\bf $CCC_{net}$} & {\bf -0.0570} & {\bf -0.1039} & {\bf 0.0032} & {\bf 0.0861} & {\bf 0.0715}
\\\hline
\end{tabular}
}
\subcaptionbox{H{\'e}non}{
\begin{tabular}{|c|C{1.5cm}|C{1.5cm}|}\hline
\diagbox{To}{From} & {\bf $x$} & {\bf $y$} \\ \hline
{\bf $x$} & 0 & -0.0273  \\ \hline
{\bf $y$} & -0.0600 & 0 \\ \hline
{\bf $CCC_{net}$} & {\bf -0.0600} & {\bf -0.0273}
\\\hline
\end{tabular}
}
\end{table}
\begin{table}[h!]
\newcommand{\xmark}{\ding{55}}
\centering
\renewcommand{\arraystretch}{1.3}
\newcolumntype{C}[1]{>{\centering\arraybackslash}m{#1}}
\caption{Indication of synchronizing variables for different dynamical systems. \checkmark  implies the variable is synchronizing and \xmark  implies the variable is not synchronizing.}
\label{table_sync_variable}
\begin{tabular}
{|c|C{1.08cm}|C{1.08cm}|C{1.08cm}|C{1.08cm}|C{1.08cm}|}\hline
\diagbox{System}{Variable} & {\bf $x$} & {\bf $y$} & {\bf $z$} & {\bf $q$} & {\bf $w$}\\ \hline
{\textit{Lorenz}} & \checkmark & \checkmark & \xmark & - & - \\ \hline
{\textit{R{\"o}ssler}} & \xmark & \checkmark & \xmark & - & -\\ \hline
{\textit{$Chen$}} & \xmark & \checkmark & \xmark & - & -
\\\hline
{\textit{5D system}} & \checkmark & \checkmark & \xmark & \xmark & \xmark
\\\hline
{\textit{H{\'e}non}} & \checkmark & \xmark & - & - & - \\ \hline
\end{tabular}
\end{table}


For the above systems, intra-system conditional $CCC$ values as well as $CCC_{net}$ values from each variable are given in Table~\ref{table_5_systems}~\footnotemark[2]. For all the systems including Lorenz, 8000 time points of time series data were taken for $CCC$ estimation after removal of 2000 transients~\footnote{While Lorenz was simulated using Euler's method, R{\"o}ssler, Chen and the 5D system were simulated using the Runge Kutta fourth order method.}. Table~\ref{table_sync_variable} indicates which variables when forced lead to complete synchronization. It can be seen from Table~\ref{table_5_systems} that the $CCC_{net}$ value from the synchronizing variable to the corresponding subsystem is always negative. In fact, in all cases (except {R{\"o}ssler}), the variable with the \emph{highest negative} $CCC_{net}$ value to the subsystem always leads to synchronization. Further, any variable with a positive $CCC_{net}$ value \emph{never} leads to synchronization. However, a variable with a negative $CCC_{net}$ value but not the highest may or may not lead to complete synchronization. For instance, in case of Lorenz and the 5D system, the $x$ variable leads to synchronization while in Chen, $x$ and in H{\'e}non, $y$ do not lead to synchronization. The exact reason(s) for this is still unclear and requires further investigation.

One of the reasons why the variable with the highest negative $CCC_{net}$ value doesn't lead to synchronization for {R{\"o}ssler} could be due to the nature of its equations. The variable $x$ which has the highest negative $CCC_{net}$ does not directly depend  on itself but on the other two variables. This is not the case for synchronizing variables in any other system that we have considered. It is possible that the high negativity of $CCC_{net}$ influence from $x$ on the subsystem $y,z$ is in fact due to $y$, the only other variable which has a negative $CCC_{net}$. Not surprisingly, it is found that $y$ leads to synchronization and not $x$.

\section{Conclusions}

In this work we provide a spatial perspective for synchronization on the basis of measuring net input causal influence to slave's non-driven subsystem. We have introduced the novel concept of \emph{causal stability} and proposed the \emph{causal stability synchronization theorem} which we have proved as a necessary and sufficient condition for synchronization in chaotic continuous-time dynamical systems. Asymptotic stability of the slave's non-driven subsystem was proven to be a necessary and sufficient condition for synchronization long back in~\cite{he1992analysis} with negative lyapunov exponents of the subsystem being its empirical condition~\cite{pecora1990synchronization}. In contrast, ergodicity of dynamical processes has allowed us to formulate an equivalent spatial condition. An empirically derived condition for causally stable subsystems has also been proposed. It involves analysis of the \emph{sign} of $CCC_{net}$ from the coupled variable to its subsystem for the master using its time series data. This is an important contribution for the control of chaos in networks where we do not know the underlying mechanism and wish to inhibit/facilitate synchronization between systems. 

Future research would involve studying several other homogeneously and heterogeneously coupled dynamical systems (both continuous and discrete-time systems) and  to analyze their $CCC_{net}$ values in order to test the proposed empirical condition. For identification of synchronizing variables, we shall explore causality measures other than CCC, that can reveal `the kind of dynamical influence', that one variable has on another. For the causal stability synchronization theorem, work will be targeted towards generalizing it for different types of coupling as well as other forms of synchronization such as phase synchronization and generalized synchronization.

\section{Acknowledgements}
We wish to dedicate this work to late Prof. Prabhakar G. Vaidya, our revered teacher/mentor, who has been a constant source of inspiration for us. Aditi Kathpalia is thankful to Manipal Academy of Higher Education for permitting this research as part of the PhD programme. The authors gratefully acknowledge the financial support of ‘Cognitive Science Research Initiative’ (CSRI-DST) Grant No. DST/CSRI/2017/54 and Tata Trusts provided for this research. 


%
\providecommand{\noopsort}[1]{}\providecommand{\singleletter}[1]{#1}%

\end{document}